\documentclass[12pt,twoside]{article}
\usepackage{hyperref}
\usepackage{geometry}
\usepackage{xcolor}
\usepackage{fancyhdr}
\usepackage{mathdots}
\usepackage{bbold}

\usepackage[labelfont=bf]{caption}
\usepackage{graphicx}
\usepackage{amsmath}
\usepackage{amssymb}
\usepackage{stmaryrd}
\usepackage{eqparbox}
\usepackage{float}
\usepackage{booktabs}
\usepackage{tabularx}
\usepackage{titling}
\usepackage{tikz-cd}
\usepackage{amsthm}
\usepackage{fancyhdr}
\usepackage{lastpage}
\usepackage{titlesec}
\usepackage[T1]{fontenc}
\titlespacing*{\section}
{0pt}{0.3ex plus 1ex minus .2ex}{0.3ex plus .2ex}
\titlespacing*{\subsection}
{0pt}{0.3ex plus 1ex minus .2ex}{0.3ex plus .2ex}
\titlespacing*{\subsubsection}
{0pt}{0.3ex plus 1ex minus .2ex}{0.3ex plus .2ex}
\makeatletter
\newcommand{\extp}{\@ifnextchar^\@extp{\@extp^{\,}}}
\def\@extp^#1{\mathop{\bigwedge\nolimits^{\!#1}}}
\makeatother

\newtheorem{Thm}{Theorem}[section]
\newtheorem{Def}{Definition}[section]

\newcommand\frontmatter{%
    \cleardoublepage
  \pagenumbering{roman}}

\newcommand\backmatter{%
  \if@openright
    \cleardoublepage
  \else
    \clearpage
  \fi
   }
\date{}

\title{\Large An Introduction to Calabi-Yau Manifolds} 
\author{\normalsize Aidan Patterson
 Waterloo, Ontario, Canada, 2021}
\usepackage{lipsum}
\renewcommand\footnotemark{}
\begin{document}
\frontmatter{}
\maketitle
\begin{figure}[h]
	\centering
	\includegraphics[width=0.3\linewidth]{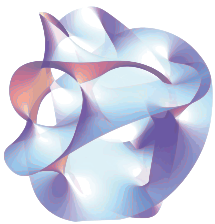}
	\caption{A cross section of the Calabi–Yau quintic in $\mathbb{CP}^4$.}
	\label{fig:springs}
\end{figure}
\begin{abstract}
  \noindent The goal of this paper is to develop the theory of Courant algebroids with integrable para-Hermitian vector bundle structures by invoking the theory of Lie bialgebroids. We consider the case where the underlying manifold has an almost para-complex structure, and use this to define a notion of para-holomorphic algebroid. We investigate connections on para-holomorphic algebroids and determine an appropriate sense in which they can be para-complex. Finally, we show through a series of examples how the theory of exact para-holomorphic algebroids with a para-complex connection is a generalization of both para-K\"{a}hler geometry and the theory of Poisson-Lie groups.
  \end{abstract}
\vspace{-9em}
\section{A Brief Introduction to Calabi-Yau Manifolds}
    \indent The first instance of the phrase ``Calabi-Yau manifolds'' occurred in 1985 in the paper ``\textit{Vacuum Configurations for Superstrings}'' by P. Candelas, G. T. Horowitz, A. Strominger, and E. Witten \cite{Vacuum Configurations for Superstrings}. In this paper, the authors place limitations on what they call the "internal manifold" $K$ of their model of space time given by $\mathcal{M}_4 \times K$, where $\mathcal{M}_4$ is the 4 dimensional Minikowski space, and $K$ is assumed to be compact. Arising from the conditions for unbroken $N=1$ supersymmetry, they deduced that this manifold must be K\"{a}hler, and additionally that it must have holonomy $SU(3)$. The reduction of the holonomy group to $SU(3)$ corresponds to the vanishing of the first Chern class (denoted by $c_1(K)$) of the manifold, and so we arrive at a definition for a Calabi-Yau manifold that is consistent with this historical development:
 \\
 \begin{Def}
 A Calabi-Yau manifold $M$ is a compact K\"{a}hler manifold such that $c_1(M)=0$.
 \end{Def}
Previous to this paper, the idea of a Calabi-Yau manifold had been conceived by Calabi and presented in a modern form during a talk at the ICM in 1954 titled: "The Space of K\"{a}hler Metrics" \cite{Proceedings}. Though he did not single out a class of manifolds for study, he was interested in the space of K\"{a}hler metrics associated with a complex manifold $M$, and how it could be parametrized. A modern form of the conjecture from this talk can be found in \cite{The Calabi Conjecture}, and is given as follows. 
\begin{Thm} (The Calabi Conjecture, 1954)
Let $M$ be a compact, complex manifold and $g$ a K\"{a}hler metric on $M$ with K\"{a}hler form $\omega$. Suppose $\tilde{R}$ is a closed, real $(1,1)$-form on $M$ with $[\tilde{R}]= c_1(M)$. Then there exists a unique K\"{a}hler metric $\tilde{g}$ on $M$ with K\"{a}hler form $\tilde{\omega}$ such that $[\tilde{\omega}]= [\omega] \in H^{2}(M,\mathbb{R})$ and the Ricci form of $\tilde{\omega}$ is $\tilde{R}$.   
\end{Thm}
Calabi managed to provide a proof of uniqueness of $\tilde{g}$, but the existence proof was completed in 1978 by Yau in \cite{Yau}. In his paper, he reduced the problem to solving an elliptic PDE of complex Monge–Amp\'{e}re type. Since this result, the field has seen some rapid development, and our goal is to cover some of the important introductory results. In particular, we will look at the interplay between the vanishing of the first Chern class and the holonomy group, and then discuss a class of examples of Calabi-Yau 3-folds that can be constructed as projective varieties in $\mathbb{CP}^N$ for some $N$ called complete intersection Calabi-Yau manifolds. Finally, we look at the special properties of the Hodge diamond for Calabi-Yau manifolds, and discuss the idea of mirror symmetry for Calabi-Yau manifolds. 
\\

\section{The Ricci Form and Chern Class}
The definition of a Calabi-Yau manifold makes reference to the first Chern class and the Ricci form of a K\"{a}hler manifold, and so we provide a brief introduction to these concepts. The curvature associated to a connection $\nabla:\Gamma(M;E)\rightarrow \Omega^1(M)\otimes \Gamma(M;E)$ on a complex $n$ dimensional vector bundle $E\rightarrow M$ is a linear map $R^{\nabla}: \Gamma(M;E)\rightarrow \Omega^2(M)\otimes \Gamma(M;E)$, given by $R^{\nabla} = \nabla \cdot \nabla$. When evaluated on sections $X,Y,Z \in \Gamma(M;E)$, one finds
\begin{align*}
R^{\nabla}(X,Y)Z &= \nabla_{X}\nabla_{Y}Z - \nabla_{Y}\nabla_{X}Z - \nabla_{[X,Y]}Z.
\end{align*}
Now, if $\{\sigma^i\}$ is a local basis for $E\rightarrow M$, then one defines the connection $1$-form $\omega$ by $\nabla \sigma^i = \omega_{j}^{i}\otimes \sigma^j$, where here we are using Einstein summation notation. In terms of $\omega$, we can compute that $R^{\nabla}(\cdot,\cdot)\sigma^i =\sum_{j}R_{ij}^{\nabla}\sigma^j =  (d\omega^{i}_j - \omega^i_k\wedge \omega^k_j)\otimes \sigma^j$. $(\Omega^{\nabla})_{j}^{i} = d\omega^{i}_j - \omega^i_k\wedge \omega^k_j$ is called the curvature 2-form.
\begin{Def}
The Ricci curvature of a connection $\nabla$ is given by $Ric(Y,Z) = Tr(X\mapsto R^\nabla (X,Y)Z)$.
\end{Def}
Now we can restrict our attention to K\"{a}hler manifolds. If $(M,g,J,\omega)$ is K\"{a}hler, then we may consider the Levi-Civita connection $\nabla^{g}$ on $M$, which since $M$ is K\"{a}hler coincides with the Chern connection (the unique connection on $M$ such that $\nabla^{(0,1)} = \bar{\partial}$). In terms of the Christoffel symbols, one finds that the non-zero entries of $Ric$ are $Ric_{i\overline{j}} = -\bar{\partial}_{j} \Gamma_{ik}^{k}$, and that $\Gamma_{ik}^{k} = \partial_{i}(\log(\sqrt{\det(g)}))$, so together, $Ric_{i\overline{j}} = -\bar{\partial}_{j}\partial_{i}(\log(\sqrt{\det(g)})) = \partial_{i}\bar{\partial}_j (\log(\sqrt{\det(g)}))$ \cite{Lectures on Kahler Geometry}. From this expression, we can define the Ricci form.
\begin{Def}
The Ricci form for a K\"{a}hler manifold $(M,g,J,\omega)$ is defined by $\rho(X,Y) = Ric(X,JY) = -i\partial \bar{\partial}\log(\sqrt{\det(g)})$.
\end{Def}
Additionally, we know that given a local basis $\{\sigma^i\}$, the Christoffel symbols can be expressed in terms of the connection $1$-form by $\Gamma_{ij}^{k} = \omega_{k}^{j}(\sigma^i)$. If we then take the trace of the curvature $2$-form $\Omega^{\nabla}$, we find $R = Tr(\Omega^{\nabla}) = \sum_{i}d\omega_{ii} = d\Gamma^{k}_{ik}$. Since $\Gamma_{ik}^{k} = \partial_{i}(\log(\sqrt{\det(g)}))$, we find that $Tr(\Omega^{\nabla}) = \bar{\partial}_j \partial_{i}\log(\sqrt{\det(g)}) = - \partial_i \bar{\partial}_j\log(\sqrt{\det(g)})$, and so $\rho = iR$. One can show that $R $ is independent of the choice of connection (see proposition 1.2 in \cite{Introduction to the Chern Class}) and since $R$ is locally $\partial \bar{\partial}$-exact, it is closed. We arrive at the following definition of the Chern class that can be used for any complex vector bundle in terms of $R$.  
\begin{Def}
Given a complex vector bundle $E\rightarrow M$, the first Chern class is $c_1(M)= \left[ \frac{i}{2\pi}R \right]\in H^{2}(M;\mathbb{R})$, where $R := Tr(\Omega^{\nabla})$ for any connection $\nabla $ on $E$. In terms of the Ricci form, we have $c_1(M)= \left[ \frac{1}{2\pi}\rho \right]$.
\end{Def}
The Calabi conjecture considers the converse of this statement, and attempts to find a K\"{a}hler structure for which $\rho \in c_1(M)$ is the Ricci form. 
\section{The Calabi Conjecture}
Given a K\"{a}hler manifold $M$, if $R,\tilde{R}\in \Omega^{(1,1)}(M)$ are in the same cohomology class, then $[\tilde{R}-R]=0$, so by the global $i\partial \bar{\partial}$-lemma, $\tilde{R} - R = \partial \bar{\partial}F$ for some smooth function $F$ on $M$. Assuming that $R$ and $\tilde{R}$ are the Ricci form corresponding to some K\"{a}hler metrics $g$ and $\tilde{g}$ respectively, components of the Ricci form are given by $R_{ij} = \partial_i \bar{\partial}_j\log(\det(g_{k\bar{l}}))$, and $\tilde{R}_{ij} = \partial_i \bar{\partial}_j\log(\det(\tilde{g}_{k\bar{l}}))$, and so one finds that
\begin{align*}
\partial \bar{\partial} \log\left( \frac{\det(\tilde{g}_{i\bar{j}})}{\det(g_{i\bar{j}})} \right) &= \partial \bar{\partial}F.
\end{align*}
Further, this tells us that $\log\left( \frac{\det(\tilde{g}_{i\bar{j}})}{\det(g_{i\bar{j}})} \right) -F$ is a constant function, and so there exists some constant $C>0$ such that 
\begin{align*}
\det(\tilde{g}_{i\overline{j}}) &= C\exp(F)\det(g_{i\overline{j}}).
\end{align*}
Finally, since we desire that K\"{a}hler forms be cohomologous, $\tilde{g}_{i\overline{j}} = g_{i\overline{j}} + \partial_{i}\bar{\partial}_j \varphi$ for some smooth function $\varphi$, and so finally, we arrive at the PDE:
\begin{align*}
\det \left( g_{i\overline{j}} + \frac{\partial^2 \varphi}{\partial z^i \partial \bar{z}^j} \right) &= C\exp(F) \det(g_{i\overline{j}}).
\end{align*}
This equation is a nonlinear, elliptic, second-order partial diﬀerential equation in $\varphi$, which means that it is of Monge-Amp\'{e}re type. This PDE is subject to the condition that 
\begin{align*}
C\int_{M}\exp(F) &= \textrm{Vol}(M).
\end{align*}
These computations give us an equivalent formulation of the Calabi Conjecture. 
\begin{Thm} (Calabi Conjecture Reformulation \cite{The Calabi Conjecture})
Let $M$ be a compact, complex manifold and $g$ a K\"{a}hler metric on $M$ with a K\"{a}hler form $\omega$. Let $F$ be a smooth function on $M$ and let $C>0$ be defined via the compatibility condition $C \int_{M}\exp(F) = \textrm{Vol}(M)$. Then there exists a unique smooth function $\varphi$ such that:
\begin{align*}
& (i) \, \, \, \, \,  \tilde{g} = \sum_{i,j} g_{i\overline{j}}  + \frac{\partial^2 \varphi}{\partial z^i  \partial \bar{z}^j} dz^i \otimes d\bar{z}^j
\\
& (ii) \, \, \, \, \, \int_{M}\varphi = 0 
\\
& (iii) \, \, \, \, \, \det\left( g_{i\overline{j}}  + \frac{\partial^2 \varphi}{\partial z^i \partial \overline{z}^j}\right) = C\exp(F)\det(g_{i\overline{j}}).
\end{align*}
\end{Thm}
A given solution $\varphi$ to the PDE above is determined up to a constant, and so the condition that $\int_{M}\varphi=0$ ensures the uniqueness of the solution. A consequence of the Calabi conjecture for Calabi-Yau manifolds is that they always admit a Ricci-flat metric.

%\begin{figure}[H]
%\includegraphics[]{figure1.png}
%\caption{\textbf{The tortoise and the hare.}}
%\end{figure}

%\begin{quotation}
 %\citep{citation1} 
 %\end{quotation}
 
\section{Holonomy of Calabi-Yau Manifolds}
We first give a brief introduction to holonomy on vector bundles. Suppose that $E\xrightarrow{\pi} M $ is a smooth vector bundle, $\nabla$ be a connection on $E$, $\{e_i\}$ a local frame, and $X \in \Gamma(M; E)$ a section locally given by $X = X^i e_i$. Given a piece-wise smooth curve $\gamma : [0,1]\rightarrow M$ with $\gamma(0)=p$ and $\gamma(1)=q$, $\sigma$ is said to be parallel along $\gamma$ with respect to $\nabla$ if $\nabla_{\dot{\gamma}(t)}\sigma = 0$ for all $t\in[0,1]$. Given some $X_p\in E_p = \pi^{-1}(p)$, one can ask if there is a section $X\in \Gamma(E;M)$ such that $\nabla_{\dot{\gamma}(t)}X = 0$, and $X\vert_{p} = X_p$. If locally $\nabla = d + \omega_{i}^{j}$, and $\gamma = (\gamma^1,\cdots, \gamma^n)$ then this is equivalent to solving the following initial value problem.
\begin{align*}
\frac{d}{dt}\left(X^j \circ \gamma\right) - \left(\omega_{i}^{j}\right)_{k}\frac{d\gamma^k}{dt} X^i\circ \gamma & = 0  & X^i\circ \gamma(0) = X_p^i 
\end{align*}
We observe that $(\omega_{i}^{j})_k$ and $\frac{d\gamma^k}{dt} $ are piece-wise smooth functions, and so the Picard-Lindel\"{o}f theorem guarantees a solution for $X_\gamma= (X^i\circ \gamma) e_i$ which can be extended to all of $[0,1]$. The parallel transport map with respect to $\gamma$ is given by $P_\gamma : E_{p}\rightarrow E_q$, $P_\gamma(X_p) = X_{\gamma(1)}$. This allows us to define the holonomy group of a manifold. 
\begin{Def}
\cite{Compact Riemannian manifolds with exceptional holonomy} Let $M$ be a connected manifold, $g$ a Riemannian metric on $TM$, and $\nabla$ the Levi-Civita connection. Suppose that $\gamma:[0,1]\rightarrow M$ is a piece-wise smooth curve with $\gamma(0)=p=\gamma(1)$. We define the Riemann holonomy group $\textrm{Hol}_p(g)$ of $g$ based at $p$ to be 
\begin{align*}
\textrm{Hol}_p(g)&:= \{P_\gamma \, : \, \gamma \textrm{ is a loop based at }p\}.
\end{align*}
$\gamma$ is called null-homotopic if it can be deformed to the constant loop at $p$. We define the restricted Riemann holonomy group:
\begin{align*}
\textrm{Hol}_p^{ \, 0}(g) &:= \{ P_\gamma \, : \, \gamma \textrm{ is a null-homotopic loop based at }p \} .
\end{align*}
\end{Def}
One immediately discovers the following facts \cite{Compact Riemannian manifolds with exceptional holonomy}. If $\gamma: [0,1]\rightarrow M$ is piece-wise smooth with $\gamma(0)=p$ and $\gamma(1)=q$, then $\textrm{Hol}_q(g) = P_\gamma \textrm{Hol}_p(g) P_\gamma^{-1}$ and $\textrm{Hol}_q^{\, 0}(g) = P_\gamma \textrm{Hol}_p^{ \, 0}(g) P_\gamma^{-1}$, and so one can speak of the holonomy group of $M$ with respect to the metric $g$ independent of the base point. Further, for $\gamma$ a loop at $p$, $P_\gamma \in O(T_pM)$, as $\nabla_{\dot{\gamma}(t)} g(P_{\gamma(t)}(v), P_{\gamma(t)}(w)) = g(\nabla_{\dot{\gamma}(t)} (P_{\gamma(t)}(v)),P_{\gamma(t)}(w)) + g(P_{\gamma(t)}(v), \nabla_{\dot{\gamma}(t)} (P_{\gamma(t)}(w))) = 0 $, as the sections are parallel, and so $g(P_{\gamma(t)}(v) , P_{\gamma(t)}(w))$ is constant along $\gamma$. This tells us that $g(v,w) = g(P_{\gamma}(v) , P_{\gamma}(w))$, as desired. 
\\

If one restricts their attention to Calabi-Yau manifolds, then the condition that $c_1(M)=0$, and hence $M$ be Ricci-flat,  places a restriction on the possible holonomy group of $M$ and vice-versa. The first restriction is as follows
\begin{Thm}
\cite{Lectures on Kahler Geometry} A compact K\"{a}hler manifold $M^{2m}$ is Calabi-Yau if and only if $\textrm{Hol}^{ \, 0}(g)\subset SU(m)$. 
\end{Thm}
Additionally, one finds restrictions that arise when the manifold is simply-connected but not necessarily compact. These are summarized in the following theorem.
\begin{Thm} \cite{Lectures on Kahler Geometry}
Let $(M,g,J,\omega)$ be a simply-connected K\"{a}hler manifold. Then the following statements are equivalent:
\begin{align*}
    &(i)  \, M \textrm{ is Ricci-flat.}
    \\
    &(ii)  \, \textrm{The canonical bundle } K=\Lambda^{(m,0)}M \textrm{ is flat}.
    \\
    &(iii)  \, \textrm{The Riemann holonomy of }M\textrm{ is a subgroup of }SU(m).
    \\
    &(iv) \, \textrm{There exists a parallel holomorphic complex volume form.}
\end{align*}
\end{Thm}
If $M$ is not simply-connected, then the final two statements are only local. Some authors take this to be the definition of a Calabi-Yau manifold, but this restriction on the full holonomy group is stronger than the assumption that $c_1(M)=0$, and is equivalent to the canonical bundle being flat. In particular, there exist manifolds with $c_1(M)=0$, but whose canonical bundles are not flat, the simplest examples of which are the hyperelliptic surfaces \cite{Compact Complex Surfaces}. The following theorem bridges the gap between theorem 4.1 and theorem 4.2 in the non-simply-connected case. 
\begin{Thm}
\cite{Lectures on Kahler Geometry} Let $M^{2m}$ be a Calabi-Yau manifold. If $m$ is odd, then $\textrm{Hol}(g) = SU(m)$, so there exists a global holomorphic $(m,0)$-form even if $M$ is not simply connected. If $m$ is even, the either $M$ is simply connected or $\pi_1(M)=\mathbb{Z}_2$ and $M$ carries no global holomorphic $(m,0)$-forms.
\end{Thm}
This result suggests that there should be a way to study general Calabi-Yau manifolds using the family of simply connected Calabi-Yau manifolds, as they often share properties. In this vein of thought, the Cheeger-Gromoll theorem \cite{Einstein Manifolds} tells us that any Calabi-Yau manifold $M$ is isomorphic to a quotient
\begin{align}\label{Cheeger-Gromoll}
M\cong (M_0 \times \mathbb{T}^l)/\Gamma,
\end{align}
where $M_0$ is a compact simply-connected K\"{a}hler manifold, $\mathbb{T}^l$ is a complex torus and $\Gamma$ is a finite group of holomorphic transformations. If one takes the De Rham decomposition of the $M_0$ from equation \ref{Cheeger-Gromoll}, we obtain $M_0 = M_1\times \cdots \times M_r$, where the $M_j$ are compact Ricci-flat simply-connected K\"{a}hler manifolds with irreducible holonomy for all $j$. Further, the Berger holonomy theorem then shows that $M_j$ has holonomy $SU(m_j)$ or $Sp(m_j)$. This gives us the following important result.
\begin{Thm} \cite{Lectures on Kahler Geometry} Any Calabi-Yau manifold $M$ is isomorphic to a quotient 
\begin{align*}
M\cong (M_1\times \cdots \times M_s \times M_{s+1}\times \cdots \times M_r \times \mathbb{T}^l)/\Gamma,
\end{align*}
where $M_j$ are compact simply-connected Calabi-Yau manifolds for $j\leq s$, and simply connected compact strict hyperk\"{a}hler manifolds for $s+1\leq j \leq r$, and $\Gamma$ is a finite group of holomorphic transformations.
\end{Thm}
At this point, it would be useful to discuss some examples of Calabi-Yau manifolds. It turns out that algebraic subsets of complex projective space are sometimes endowed with the structure of a Calabi-Yau manifold, as we will now see.
\section{Complete Intersection Calabi-Yau Manifolds}
To begin, a compact complex manifold $M$ is called projective it if can be holomorphically embedded in some complex projective space $\mathbb{CP}^N$. All projective manifolds are algebraic (by a well known result of Chow), meaning they are the intersection of the projective zero sets of a finite set of homogeneous polynomials. In the case of Calabi-Yau manifolds, we have the following remarkable result.
\begin{Thm}\cite{Lectures on Calabi-Yau and Special Lagrangian Geometry}
Every Calabi-Yau manifold of complex dimension $m\geq 3$ is projective. That is, given a K\"{a}hler manifold $(M,g,J,\omega)$, the complex manifold $(M,J)$ is biholomorphic to a complex projective variety of $\mathbb{CP}^N$.
\end{Thm}
One way of constructing Calabi-Yau manifolds is, therefore, to look at complex algebraic varieties in $\mathbb{CP}^N$, and check against topological invariants to see which of these varieties are Calabi-Yau. This will not necessarily produce all Calabi-Yau manifolds, but can produce a Calabi-Yau manifold with the same underlying complex manifold as any Calabi-Yau manifold. A rather simple set of conditions on the degrees of the homogeneous polynomials defining a complex projective variety endows the variety with a Calabi-Yau structure.
\begin{Thm} \cite{Lectures on Calabi-Yau and Special Lagrangian Geometry}
Let $P_1,\cdots, P_k$ be generic homogeneous irreducible polynomials in $m+1$ variables, and $d_1,\cdots, d_k$ represent the degrees of these polynomials. Then $N= \{[z_0:\cdots:z_{m+1}]\in \mathbb{CP}^{m}: P_i=0\}$ is Calabi-Yau if and only if  $m\geq 3+k$ and $\sum_i d_i =m+1$. In this case, we call $N$ a complete intersection Calabi-Yau (CICY) manifold, and it has dimension $m-k$.
\end{Thm}
Now, if we at first consider degree $1$ polynomials, it can be shown that their zero set is isomorphic to $\mathbb{CP}^{N-1}$, and so one only needs to consider the degree $2$ or higher curves to come up with new examples. If we wish to enumerate possible classes of CICY 3-folds, then we simply find lists of $k=m-3$ positive integers greater than $1$ that sum up to $m+1$. We denote these spaces as $Y_{m;d_1,\cdots,d_k}$, and can construct the following exhaustive list \cite{Superstring Theory}. 
\begin{align}\label{Table}
    \begin{array}{c|c|c}
      \textrm{Ambient Space}   & \textrm{Polynomial Degrees} & \textrm{Euler Characteristic} \\
    \hline  \mathbb{CP}^4   & (5) &\chi (Y_{4;5})=-200 \\
    \mathbb{CP}^5 & (4,2), \,  \, (3,3) & \chi(Y_{5;4,2})=-176, \, \,  \, \chi(Y_{5;3,3}) = -144 \\
    \mathbb{CP}^{6} & (3,3,2) & \chi(Y_{6;3,3,2}) = - 144\\
    \mathbb{CP}^{7} & (2,2,2,2) & \chi(Y_{7;2,2,2,2})= -128
    \end{array}
\end{align}
We have included the Euler characteristics for later use. In general, this process will always produce a finite number of topological types for each dimension. There are ways to construct additional examples in a similar way by considering the ambient space $\mathbb{CP}^{N_1}\times \cdots \times \mathbb{CP}^{N_F}$, and considering the zero sets of polynomials which are homogeneous within each $\mathbb{CP}^{N_j}$ factor, which may still be referred to as CICY manifolds. In \cite{Complete Intersection Calabi-Yau Manifolds}, the authors constructed matrices out of the lists of degrees of polynomials where each row corresponds to the degrees of homogeneous polynomials in each factor, and used topological obstructions to show that there were 7868 possible pairs of ambient spaces and general homogeneous polynomials that produce Calabi-Yau 3-folds. Further, it was shown in \cite{All the Hodge Numbers} that there were at least 265 distinct Hodge diamonds in this collection, and so at least 265 topologically distinct CICY 3-folds that can be constructed this way. We now move on to discussing the Hodge diamond for Calabi-Yau 3-folds. 
\section{Hodge Diamonds of Calabi-Yau 3-folds and Mirror Symmetry}
To begin, we recall that given a complex manifold $M$, there exists a decomposition of the differential forms on $M$ into forms of type $(p,q)$ with $\Omega^k (M) = \sum_{p+q=k}\Omega^{(p,q)}(M)$, and a differential operator $\bar{\partial} : \Omega^{(p,q)}(M)\rightarrow \Omega^{(p,q+1)}(M)$ with $\bar{\partial}^2=0$. We define the Betti numbers of $M$ to be $b_k = \dim \left( H^{k}_{dR}(M) \right)$), and using the fact that $\bar{\partial}^2 =0$, we define the Hodge numbers to be $h^{p,q} = \dim \left( H^{(p,q)}_{\bar{\partial}}(M)  \right)$. We have the following result for compact K\"{a}hler manifolds. 
\begin{Thm}
\cite{Lectures on Kahler Geometry} The following relations hold for the Betti and Hodge numbers on a compact K\"{a}hler manifold $M^{2m}$ for all $1\leq k \leq 2m$, and $1\leq p,q\leq m$.
\begin{align*}
b_k&=b_{2m-k} & b_k&= \sum_{p+q=k}h^{p} & h^{p,q} &= h^{q,p} & h^{p,q}&= h^{m-p,m-q} & h^{p,p}\geq 1 
\end{align*}
In particular, this shows that all the Betti numbers of odd order are even, and all the Betti numbers of even order are non-zero. 
\end{Thm}
The Hodge diamond is obtained by writing the Hodge numbers in rows whose sum is a Betti number. In general, this takes the form
\begin{align*}
\begin{array}{ccccccccc}
& & & & h^{m,m} & & & &   \\
& & & h^{m,m-1} &  &h^{m-1,m} & & & \\
& &h^{m,m-2} & &h^{m-1,m-1} & & h^{m-2,m} & &  \\
&\iddots & \vdots & &\vdots & & \vdots&\ddots &  \\
h^{m,0}& &h^{m-1,1} &\cdots & & \cdots& h^{1,m-1}& & h^{0,m} \\
&\ddots & \vdots& &\vdots & & \vdots &\iddots &  \\
& &h^{2,0} & & h^{1,1} & & h^{0,2}& &  \\ 
& & & h^{1,0}& &h^{0,1} & & &  \\
& & & & h^{0,0} & & & &  \\
\end{array}
\end{align*}
We are particularly interested in the possible Hodge diamonds of a CICY 3-fold. In this case, the fact that these manifolds have 3 complex dimensions implies that the canonical bundle of each of these manifolds is trivial (by theorem 4.3), and so there exists a unique (up to scaling) holomorphic volume form $\Omega$. In terms of the Hodge numbers, this implies that $h^{3,0}=1$. Further given a $(0,q)$ cohomology class $[\alpha]$, there is a unique $(0,3-q)$ cohomology class $[\beta]$ such that $\int_M \alpha \wedge \beta \wedge \Omega =1$  (see \cite{Lectures on Complex Geometry}). Thus $h^{0,q}=h^{0,3-q}$. Therefore, $h^{3,0}=h^{0,3} = h^{0,0} = h^{3,3}=1$. One can also show that $h^{1,0}=0$ \cite{Beastiary}, and so $h^{1,0}=h^{0,1} = h^{0,2} = h^{2,0}=h^{2,3} = h^{3,2} = h^{3,1} = h^{1,3} = 0 $. Therefore the only independent Hodge numbers are $h^{1,1} =h^{2,2}$ and $h^{2,1}=h^{1,2}$. The Hodge diamond is given by 
\begin{align*}
\begin{array}{ccccccc}
     & & &1 & & &  \\
     & & 0& & 0& &  \\
     & 0& &h^{1,1} & & 0&  \\
    1 & & h^{2,1}& & h^{2,1}& & 1 \\
     &0 & &h^{1,1} & &0 &  \\
     & &0 & & 0& & \\
     & & &1 & & & \\
\end{array}
\end{align*}
Now, the last two unspecified Hodge numbers are related via the Euler characteristic. The Euler characteristic is given as the alternating sum of the Betti numbers, and so $\chi(M)= \sum_{k=1}^{2m}(-1)^k b_k = 2b_0 - 2b_1 + 2b_2 - b_3= 2(h^{1,1}-h^{2,1})$. More concretely, $h^{1,1}$ classifies the infinitesimal deformations of the K\"{a}hler structure (representing a choice of the Ricci-form), and $h^{2,1}$ represents the possible deformations of the complex structure (Chapter 6 \cite{Mirror Symmetry}). In order to compute these invariants for the case of CICY 3-folds that are subsets of a single $\mathbb{CP}^N$, as in the Table \ref{Table}, one uses the Lefschetz Hyperplane theorem (\cite{Beastiary}, Page 45), which in this case states that $H^{(p,q)}(\mathcal{M},\mathbb{Z}) = H^{(p,q)}(\mathbb{CP}^N,\mathbb{Z})$ for $p+q\neq 3$. This allows us to compute $h^{1,1}$ using the Hodge numbers for the ambient space, which in all of these cases reduces to $h^{1,1}=1$. Famously, this says that for the CICY quintics in $\mathbb{CP}^4$, $h^{2,1} = 101$. We obtain the Hodge diamond 
\begin{align*}
&Y_{4;5}: & \begin{array}{ccccccc}
     & & &1 & & &  \\
     & & 0& & 0& &  \\
     & 0& &1 & & 0&  \\
    1 & & 101& & 101& & 1 \\
     &0 & &1 & &0 &  \\
     & &0 & & 0& & \\
     & & &1 & & & \\
\end{array}&
\end{align*}
Once again, we take a look at the invariants $h^{1,1}$ and $h^{2,1}$. The principle of mirror symmetry for Calabi-Yau 3-folds is this: Given a Calabi-Yau 3-fold $M$, there is often a way to associate to it a "mirror" Calabi-Yau 3-fold $\tilde{M}$ with $H^{(1,1)}(M) = H^{(2,1)}(\tilde{M})$ and $H^{(2,1)}(M) = H^{(1,1)}(\tilde{M})$. For the given example of the CICY quintic $ \sum_i z_i^5 = 0$, one would expect there to be a Calabi-Yau 3-fold $\tilde{Y}_{4;5}$ with $h^{1,1} = 101$ and $h^{2,1}= 1$. This manifold in particular is constructed in \cite{Mirror Symmetry} as follows. They consider a 1-parameter family of quintics: $\sum_{i}X_i^5 -5\psi \prod_{i}X_i$, and note that this is preserved under the action of a fifth root of unity $\lambda$ by $X_i \mapsto \lambda^{k_i}X_i$ with $\sum_{i}k_i \equiv 0 \mod 5$. In fact, when one considers that $\mathbb{CP}^4$ is invariant under scaling, this can be augmented to an action of $(\mathbb{Z}_5)^3$. The manifold $\tilde{Y}_{4;5}$ is then defined as 
\begin{align*}
\tilde{Y}_{4;5} = \left. \left( \sum_{i}X_i^5 -5\psi \prod_{i}X_i \right)\middle/ (\mathbb{Z}_5)^3 \right. .
\end{align*}
The action of $(\mathbb{Z}_5)^3$ has fixed points however, and so this manifold is singular. However, there exists a method in this case to resolve those singularities, and when all is said and done, the final smooth $\tilde{Y}_{4;5}$ has the desired Hodge numbers. In general, it's expected that for each Calabi-Yau manifold $M$ with Hodge numbers $h^{p,q}$, there exists a Calabi-Yau manifold $\tilde{M}$ with Hodge numbers $\tilde{h}^{p,q} =  h^{m-p,q}$. The question of whether this can be done in general is still open, and is an active area of research.
\pagebreak

\end{document}